\documentclass[12pt]{amsart}
\usepackage{latexsym,amssymb,amsmath}
\usepackage{amsmath,amsfonts}
\usepackage{epsfig}
\usepackage{graphicx}
\usepackage{verbatim}
\newtheorem{proposition}{Proposition}
\newtheorem{theorem}{Theorem}
\newtheorem{lemma}{Lemma}

\allowdisplaybreaks

\def\nn{\nonumber}
\def\T{\mathbb T}
\def\mn{\mathcal N}

\def\R{{\mathbb R}}
\def\Z{{\mathbb Z}}
\def\la{\langle}
\def\ra{\rangle}

\def\les{\lesssim}
\def\1{{\bf 1}}

\def\eqnn{\begin{eqnarray*}}
\def\eeqnn{\end{eqnarray*}}
\def\eqn{\begin{eqnarray}}
\def\eeqn{\end{eqnarray}}
\newcommand{\nc}{\newcommand}
\nc{\be}{\begin{equation}}
\nc{\ee}{\end{equation}}
\nc{\ba}{\begin{eqnarray}}
\nc{\ea}{\end{eqnarray}}
\nc{\eps}{\epsilon}
\def\prf{\begin{proof}}
\def\endprf{\end{proof}}
\begin{document}

\title{Talbot effect for the cubic nonlinear Schr\"odinger equation on the torus}

\author{{\bf M.~B.~Erdo\u gan, N.~Tzirakis}\\
University of Illinois\\
Urbana-Champaign}

\thanks{The authors were partially supported by NSF grants DMS-1201872 (B.~E.), and DMS-0901222 (N.~T.) }

\date{}

\begin{abstract}
We study the evolution of the one dimensional periodic cubic Schr\"odinger equation (NLS) with bounded variation data. For the linear evolution, it is  known   that for irrational times the solution is a continuous, nowhere differentiable fractal-like curve. For rational times the solution is a linear combination of finitely many translates of the initial data. Such a dichotomy was fist observed by Talbot in an optical experiment  performed in 1836, \cite{talbot}. In this paper we prove that a similar phenomenon occurs in the case of the   NLS equation.
\end{abstract}

\maketitle
\section{Introduction}

In a 1836 optical experiment  Talbot \cite{talbot}  observed white light passing through a diffraction grating. He looked at the images that were produced with the help of magnifying lens and noticed a sharp focused pattern with a certain periodicity depending on the distance.
Rayleigh  \cite{ray}  calculated the Talbot distance as $d=\frac{\alpha^2}{\lambda}$ where $\alpha$ is the spacing of the grating and $\lambda$ is the wavelength of the incoming light.

Berry with his collaborators (see, e.g., \cite{mber, berklei,berlew, bermar})   studied the Talbot effect in a series of papers.
In particular, in \cite{berklei},  Berry and Klein  used the linear Schr\"odinger evolution to model the Talbot effect. They showed that at rational multiples of the Talbot distance overlapping copies of the grating pattern reappear, while at irrational times the images have a fractal nowhere differentiable profile. Also in  \cite{mber}, Berry  conjectured that  for the $n-$dimensional linear Schr\"odinger equation confined in a box
 the imaginary part $\Im u(x,t)$, the real part $\Re{u(x,t)}$ and the density $|u(x,t)|^2$ of the solution is a fractal set with dimension $D=n+\frac{1}{2}$ for most irrational times.  He also observed that in the one dimensional case there are space slices whose time fractal dimension is $\frac74$ and there are diagonal slices with dimension $\frac54$.

The first mathematically rigorous work in this area appears to be due to Oskolkov. In \cite{osk}, he studied  a large class of linear dispersive equations with bounded variation initial data.    In the case of the linear Schr\"odinger equation, he proved that  at rational times the solution is a continuous function  of $x$ and at rational times it is a bounded function with at most countably many discontinuities. For the exact statement see Theorem \ref{oskolkov} below.

The idea that the profile of linear dispersive equations depend on the algebraic properties of time have been further exploited in the papers of Kapitanski-Rodniaski  \cite{rodkap}, Rodnianski  \cite{rod}, and Taylor  \cite{mtay}. In \cite{rodkap}, the authors show that the solution  to the linear Schr\"odinger equation has better regularity properties (measured in Besov spaces) at irrational than rational times. It is important to note that this subtle effect can not be observed  in the scale of Sobolev spaces  since the linear propagator is a unitary operator in  Sobolev spaces. In \cite{rod}, using the result in \cite{rodkap}, Rodnianski partially justified  Berry's conjecture in one dimension, see Section~\ref{sec:discuss}.

In \cite{mtay}, Taylor independently obtained Berry's quantization result and also extended it to higher dimensional spheres and  tori \cite{mtay2}.
In particular he  proved that at rational times the solution is a linear combination of finitely many translates of the initial data with the coefficients being Gauss sums. He further showed that some classical identities for Gauss sums can be obtained by an analysis of the linear Schr\"odinger evolution.
This shouldn't come as a surprise since number theoretic connections of Talbot effect had already been observed in the works of Oskolkov, Kapitanski and Rodnianski. In \cite{mtay},
Taylor also noted that the quantization implies the $L^p$ boundedness of the multiplier $e^{it\Delta}$ for rational values of $\frac{t}{2\pi}$. It is  known that, \cite{mtay}, the propagator is unbounded in $L^p$ for $p\neq 2$ and $\frac{t}{2\pi}$ irrational.  This can be considered as another manifestation of the Talbot effect.

More recently   Olver \cite{olv}, and  Chen and Olver  \cite{chenolv,chenolv1}  provided numerical simulations of the Talbot effect for a large class of dispersive equations, both linear and nonlinear. In the case of polynomial dispersion, they numerically confirmed the rational/irrational dichotomy discussed above. This behavior persists for both integrable and nonintegrable systems. An important question that the authors raised is the appearance of such phenomena in the case of nonpolynomial dispersion relations. The numerics demonstrate that the large wave number asymptotics of the dispersion relation plays the dominant role governing the qualitative features of the solutions.  We should also note that in \cite{ZWZX} the Talbot effect was observed experimentally in a nonlinear setting.

In this note we investigate the Talbot effect for  cubic  nonlinear  Schr\"odinger equation (NLS) with periodic boundary conditions.  Our goal is to extend Oskolkov's and Rodnianski's results for bounded variation data to the NLS evolution, and  provide rigorous confirmation of some numerical observations in \cite{olv,chenolv,chenolv1}.
In particular we prove that for a large class of rough data the solution of the  NLS equation, for almost all times, is a continuous but fractal-like curve with upper Minkowski dimension\footnote{Upper Minkowski (also known as fractal) dimension, $\overline{\text{dim}}(E)$, of a bounded set $E$ is given by $$\limsup_{\epsilon\to 0}\frac{\log({\mathcal N}(E,\epsilon))}{\log(\frac1\epsilon)},
$$
where ${\mathcal N}(E,\epsilon)$ is the minimum number of $\epsilon$--balls required to cover $E$.}
$\frac{3}{2}$.
Our main result is the following theorem:
\begin{theorem}\label{mainthm}
Consider the nonlinear Schr\"odinger equation on the torus:
\begin{align*}
&iu_t+u_{xx}+|u|^2u=0,\,\,\,\,\,\, t\in \R,\,\,\,x\in\T=\R/2\pi\Z,\\
&u(x,0)=g(x).
\end{align*}
Assuming that $g$ is of bounded variation, we have\\
i)  $u(x,t)$ is a continuous function of $x$ if $\frac{t}{2\pi}$ is an irrational number. For rational values of $\frac{t}{2\pi}$, the solution is a bounded function with at most countably many discontinuities.  Moreover, if $g$ is also continuous  then $u\in C^0_tC^0_x$.\\
ii) If in addition $g\not \in \bigcup_{\epsilon>0}H^{\frac12+\epsilon}$,
then for almost all times    either the real part or the imaginary part of the graph of
$u(\cdot,t)$ has upper Minkowski dimension $ \frac32$.
\end{theorem}

We note that the simulations in \cite{olv,chenolv,chenolv1} were performed in the case when $g$ is a step function, and that Theorem~\ref{mainthm} applies in that particular case.

To prove Theorem~\ref{mainthm} we first obtain a smoothing result for NLS stating that the nonlinear Duhamel part of the evolution is smoother than the linear part by almost half a derivative. For bounded variation data, this implies that the nonlinear part is in $H^{1-}$  which immediately yields the upper bound on the dimension of the curve.  The lower bound is obtained by combining our smoothing estimate with Rodnianski's result in \cite{rod},
and an observation from \cite{DJ} connecting smoothness and geometric dimension. We remark that the first part of
Theorem~\ref{mainthm} was observed in \cite{erdtzi} in the case of KdV equation.

\section{Discussion of earlier results and the proof of Theorem~\ref{mainthm}}  \label{sec:discuss}
First of all recall that for $s\geq 0$, $H^s(\T)$ is defined as a subspace of $L^2$ via the norm
$$
\|f\|_{H^s(\T)}:=\sqrt{\sum_{k\in\Z} \la k\ra^{2s} |\widehat{f}(k)|^2},
$$
where $\la k\ra:=(1+k^2)^{1/2}$ and $\widehat{f}(k)=\int_0^{2\pi}f(x)e^{-ikx} dx$ are the Fourier coefficients of $f$. We note that the local and global wellposedness of NLS for $H^s$ data for $s\geq 0$ was obtained by Bourgain \cite{bourgain}.

We start by formally decomposing the solution as
$$
u(x,t)=e^{i(\partial_{xx}+P)t} g+\mn(x,t),
$$
where $P=\|g\|_2^2/\pi$. Here $\mn$ is the nonresonant part of the nonlinear Duhamel term of the solution.
The proof of Theorem~\ref{mainthm} relies on earlier results by Oskolkov and by Rodnianski on the linear part $e^{i(\partial_{xx}+P)t} g$ and the
following smoothing theorem involving $\mn$.
\begin{theorem}\label{mainpro}
Fix $s>0$. Assume that $g\in H^s(\T)$. Then for any $a<\min(2s,1/2)$,  we have
$$\mn(x,t) \in C^0_{t\in\R}H^{s+a}_{x\in\T}.$$
\end{theorem}

\noindent
We note that this type of smoothing does not hold for $s\leq 0$. For this and (local) smoothing results in $\mathcal F \ell^p$ spaces, see  \cite{chr}.
We will prove Theorem~\ref{mainpro} in Section~\ref{sec:t2} below.
Note that if $g$ is of bounded variation then $g\in\bigcap_{\epsilon>0}H^{\frac12-\epsilon}$, and hence
$$\mn(x,t)  \in \bigcap_{\epsilon>0} C^0_{t\in\R}H^{1-\epsilon}_{x\in\T}.$$
In particular,  $\mn(x,t)  \in C^0_tC^0_x$.
This and the following theorem\footnote{In fact Theorem~\ref{oskolkov}  is a special case of a Theorem of Oskolkov which asserts the statement above for one dimensional dispersive equations with polynomial dispersion relation.} of Oskolkov  \cite{osk} conclude  the proof of part i) of Theorem~\ref{mainthm}.

\begin{theorem} \cite{osk} \label{oskolkov}
Let   $g$ be of bounded variation, then $e^{it\partial_{xx}}g$ is a continuous function of $x$ if $\frac{t}{2\pi}$ is an irrational number. For rational values of $\frac{t}{2\pi}$, it is a bounded function with at most countably many discontinuities.  Moreover, if $g$ is also continuous  then $e^{it\partial_{xx}}g\in C^0_tC^0_x$.
\end{theorem}

\vskip 0.1in
\noindent
{\bf Remark 1.} In the case when the initial data has jump discontinuities, by the quantization result in
\cite{berklei, mtay, olv}, the linear solution has jump discontinuities when $\frac{t}{2\pi}$ is a rational number, since it is a linear combination of finitely many translates of the initial data. Noting that  $\mn$ is a continuous function of $x$, we have the same conclusion for the NLS evolution.

The second part of the theorem will rely on a result of Rodnianski \cite{rod} which in turn relies on results obtained by Kapitanski  and Rodnianski \cite{rodkap}. To state this result we need to define the Besov space $B^s_{p,\infty}$ via the norm:
$$
\|f\|_{B^s_{p,\infty}}:=\sup_{j\geq 0}2^{sj} \|P_j f\|_{L^p},
$$
where $P_j$ is a Littlewood-Paley projection on to the frequencies $\approx 2^j$. We should note that $C^\alpha(\T)$ coincides with $B^\alpha_{\infty,\infty}$, see, e.g., \cite{tri}.

\begin{theorem}\cite{rod}\label{rod}
Assume that $g$ is of bounded variation and that $g\not \in \bigcup_{\epsilon>0}H^{\frac12+\epsilon}$, then for almost all irrational $\frac{t}{2\pi}$,
\begin{align*}
&e^{it\partial_{xx}}g\not\in\bigcup_{\epsilon>0} B^{\frac12+\epsilon}_{1,\infty},\\
&e^{it\partial_{xx}}g\in\bigcap_{\epsilon>0} B^{\frac12-\epsilon}_{\infty,\infty}=\bigcap_{\epsilon>0} C^{\frac12-\epsilon}.
\end{align*}
\end{theorem}

We  note that for each $\alpha$, $H^\alpha\in B^\alpha_{1,\infty}$, and that, for $0<\alpha<\frac12$,  $H^{\alpha+\frac12}\subset C^\alpha$. Therefore by Theorem~\ref{mainthm}, we have for each $t$,
\begin{align*}
&\mn(\cdot,t) \in\bigcap_{\epsilon>0} B^{1-\epsilon}_{1,\infty},\\
&\mn(\cdot,t)\in \bigcap_{\epsilon>0} C^{\frac12-\epsilon}.
\end{align*}
This implies that Theorem~\ref{rod} is valid for the  nonlinear solution $u$.

It is a well-known result that if $f:\T\to\R$ is in $C^\alpha$, then
the graph of $f$ has upper Minkowski dimension $D\leq 2-\alpha$. Therefore, the graphs of $\Re(u)$ and $\Im(u)$ have dimension at most $\frac32$. We note that this upper bound wouldn't have followed if we had less than half  a  derivative gain in Theorem~\ref{mainpro}.

The lower bound follows from the claim of Theorem~\ref{rod} for the nonlinear solution $u$ and the following theorem of Deliu and Jawerth \cite{DJ}.

\begin{theorem}\cite{DJ} The graph of a continuous function $f:\T\to\R$  has upper Minkowski dimension $D \geq 2-s$ provided that $f\not\in\bigcup_{\epsilon>0} B^{s+\epsilon}_{1,\infty}$.
\end{theorem}

\section{Proof of Theorem~\ref{mainpro}}\label{sec:t2}

The theorem follows from a simple $X^{s,b}$ space estimate that appears to be new.

With the change of variable $u(x,t)\to u(x,t) e^{iPt}$, where $P=\|g\|_2^2/\pi$, we obtain the equation
$$
iu_t+u_{xx}+|u|^2u-Pu=0,\,\,\,\,\,\, t\in \R,\,\,\,x\in\T,
$$
with initial data in $g\in H^s$, $s>0$.
We want to prove that  for any $a<\min(2s,\frac{1}{2})$,  we have
$$u-e^{i\partial_{xx}t}g \in C^0_{t\in\R}H^{s+a}_{x\in\T}.$$

Note the following identity which follows from Plancherel's theorem:
\begin{multline*}
\widehat{|u|^2u}(k)=\sum_{k_1,k_2 }   \widehat u(k_1)\overline{\widehat u(k_2)}\widehat u(k-k_1+k_2)\\ =\frac1\pi
\|u\|_2^2\widehat{u}(k)-|\widehat u(k)|^2\widehat{u}(k)+\sum_{k_1\neq k, k_2\neq k_1}  \widehat u(k_1)\overline{\widehat u(k_2)}\widehat u(k-k_1+k_2)\\
=: P  \widehat{u}(k) +\widehat{\rho(u)}(k)+ \widehat{R(u)}(k),
\end{multline*}
where $\widehat{u}(k)=\int_0^{2\pi}u(x)e^{-ikx} dx$ are the Fourier coefficients.
Using this in the Duhamel's formula, we have
$$
u(t)=e^{it \partial_{xx} }g+i\int_0^t e^{i(t-\tau) \partial_{xx} }(\rho(u)+R(u))d\tau.
$$
We note that
\be\label{rhobound}
\|\rho(u)\|_{H^{s+a}}=\sqrt{\sum_k |\widehat{u}(k)|^6 \la k\ra^{2s+2a}}\les \|u\|_{H^s}^3,
\ee
for $0\leq a\leq 2s$.

Using \eqref{rhobound}, we have
$$
\|u(t)-e^{it \partial_{xx} }g\|_{H^{s+a}}\les \int_0^t \|u(\tau)\|_{H^s}^3 d\tau+\Big\|\int_0^t e^{i(t-\tau) \partial_{xx} } R(u) d\tau\Big\|_{H^{s+a}}.
$$
To bound the integral involving $R(u)$, we work with the $X^{s,b}$ space \cite{bourgain,Bou2}:
$$
\|u\|_{X^{s,b}}=\big\| \widehat{u}(\tau,k)\la k\ra^{s} \la \tau-k^2\ra^{b} \big\|_{L^2_\tau\ell^2_k},
$$
with $b>1/2$. We also define the restricted norm
$$
\|u\|_{X^{s,b}_\delta}=\inf_{ \tilde u=u \text{ on } [-\delta,\delta]}\|\tilde u\|_{X^{s,b}}.
$$
We will use the embedding $X^{s,b}\subset C^0_t H^{s}_x$ for $b>\frac{1}{2}$ and the following inequality from \cite{bourgain}.
For any $s\in \mathbb R$, $\delta \leq 1$, and $b>\frac{1}{2}$, we have
$$\Big\| \int_0^t e^{-i(t-s)\partial_{xx}} F(s) ds \Big\|_{X^{s,b}_{\delta}}\lesssim \|F\|_{X^{s,b-1}_\delta}.$$

\noindent
Thus for $0\leq t\leq \delta$, we have
$$
\Big\|\int_0^t e^{i(t-\tau) \partial_{xx} } R(u) d\tau\Big\|_{H^{s+a}}\les \Big\|\int_0^t e^{i(t-\tau) \partial_{xx} } R(u) d\tau\Big\|_{X^{s+a,b}_\delta}
\les \|R(u)\|_{X^{s+a,b-1}_\delta}.
$$
We will prove the following proposition in the next section.
\begin{proposition}\label{lem:smooth} For fixed $s>0$ and $a<\min(2s,\frac{1}{2})$, we have
$$\|R(u)\|_{X^{s+a,b-1}}\les \|u\|_{X^{s,b}}^3.$$
provided that $0<b-\frac{1}{2}$ is sufficiently small. The same inequality holds for the restricted norms.
\end{proposition}

Using Proposition~\ref{lem:smooth}, we have for $0<t<\delta<1$ (where $[0,\delta]$ is the local existence interval)
\be\label{localsm}
\|u(t)-e^{it \partial_{xx} }g\|_{H^{s+a}}\les \int_0^t \| u(\tau) \|_{H^s}^3d\tau+\|u\|_{X^{s,b}_\delta}^3\les \|u\|_{X^{s,b}_\delta}^3 \les \|g\|_{H^s}^3.
\ee
We note that (see page 113-114 in \cite{Bbook}) the local existence time $\delta$ depends on the $L^2$ norm of the data, $\delta=\delta(\|g\|_2)<1$, and by iteration we also have for any $s\geq 0$,
$$
\|u\|_{H^s}\leq C e^{C|t|} \|g\|_{H^s}=:T(t).
$$

Fix $t$ large.  For $r\leq t$, we have the bound
$$\|u(r)\|_{H^s}\lesssim T(r)\leq T(t).$$
Using \eqref{localsm} repeatedly, we have
$$
\|u(j\delta)-e^{i\delta \partial_{xx}}u((j-1)\delta)\|_{H^{s+a}}\les \|u((j-1)\delta)\|_{H^s}^3 \lesssim T(t)^{3},
$$
for any $j\leq t/\delta$.
Using this we obtain (with $J=t/\delta$)
\begin{align*}
\|u(J\delta)-e^{iJ \delta\partial_{xx}}f\|_{H^{s+a}}
&\leq \sum_{j=1}^J\|e^{i(J-j) \delta \partial_{xx} }u(j\delta)-e^{i(J-j+1) \delta \partial_{xx}}u((j-1)\delta)\|_{H^{s+a}}\\
&= \sum_{j=1}^J\| u(j\delta)-e^{ i\delta \partial_{xx}}u((j-1)\delta)\|_{H^{s+a}}\lesssim J T(t)^3 \approx t T(t)^3/\delta.
\end{align*}
This finishes the proof of Theorem~\ref{mainpro}.

\section{Proof of Proposition~\ref{lem:smooth}}
We start with the following elementary lemma (see, e.g.,  the Appendix of \cite{erdtzi1}).
\begin{lemma}\label{lem:sums}   If  $\beta\geq \gamma\geq 0$ and $\beta+\gamma>1$, then
\be\nn
\sum_n\frac{1}{\la n-k_1\ra^\beta \la n-k_2\ra^\gamma}\lesssim \la k_1-k_2\ra^{-\gamma} \phi_\beta(k_1-k_2),
\ee
where
 \be\nn
\phi_\beta(k):=\sum_{|n|\leq |k|}\frac1{|n|^\beta}\sim \left\{\begin{array}{ll}
1, & \beta>1,\\
\log(1+\la k\ra), &\beta=1,\\
\la k \ra^{1-\beta}, & \beta<1.
 \end{array}\right.
\ee
\end{lemma}

Using the definition of $X^{s,b}$ norm, we have
$$
\|R(u)\|_{X^{s+a,b-1}}^2 =\Big\|\int_{\tau_1,\tau_2}\sum_{k_1\neq k,k_2\neq k_1 } \frac{\la k \ra^{s+a} \widehat u(k_1,\tau_1)\overline{\widehat u(k_2,\tau_2)}\widehat u(k-k_1+k_2,\tau-\tau_1+\tau_2)}{ \la \tau-k^2\ra^{1-b}}\Big\|_{\ell^2_kL^2_\tau}^2.
$$
Let
$$
f(k,\tau)=|\widehat u(k,\tau)| \la k\ra^{s } \la \tau-k^2\ra^{b}.
$$
It suffices to prove that
\begin{multline*}
\Big\|\int_{\tau_1,\tau_2}\sum_{k_1\neq k,k_2\neq k_1 } M(k_1,k_2,k,\tau_1,\tau_2,\tau) f(k_1,\tau_1)f(k_2,\tau_2)f(k-k_1+k_2,\tau-\tau_1+\tau_2) \Big\|_{\ell^2_kL^2_\tau}^2\\ \lesssim \|f\|_2^6=\|u\|_{X^{s,b}}^6,
\end{multline*}
where
\begin{multline}\nn
M(k_1,k_2,k,\tau_1,\tau_2,\tau)=\\
\frac{\la k \ra^{s+a}\la k_1\ra^{-s } \la k_2\ra^{-s } \la k-k_1+k_2\ra^{-s }  }{  \la \tau-  k^2\ra^{1-b}\la \tau_1- k_1^2\ra^{b}\la \tau_2-k_2^2\ra^{b}\la \tau-\tau_1+\tau_2- (k-k_1+k_2)^2\ra^b}.
\end{multline}
By Cauchy--Schwarz in $\tau_1,\tau_2,k_1,k_2$ variables, we estimate the norm above by
\begin{multline*}
\sup_{k,\tau}\Big( \int_{\tau_1,\tau_2}\sum_{k_1\neq k,k_2\neq k_1 }  M^2(k_1,k_2,k,\tau_1,\tau_2,\tau) \Big)\times \\ \Big\|\int_{\tau_1,\tau_2}\sum_{k_1,k_2 }  f^2(k_1,\tau_1)f^2(k_2,\tau_2)f^2(k-k_1+k_2,\tau-\tau_1+\tau_2) \Big\|_{\ell^1_kL^1_\tau}.
\end{multline*}
Note that the norm above is equal to
$\big\|  f^2*f^2*f^2 \big\|_{\ell^1_kL^1_\tau}$, which can be estimated by $ \|f\|_2^6 $
by Young's inequality. Therefore, it suffices to prove that the supremum above is finite.

Using  Lemma~\ref{lem:sums} in $\tau_1$ and $\tau_2$ integrals, we obtain
\begin{align*}
\sup_{k,\tau} \int_{\tau_1,\tau_2}\sum_{k_1\neq k,k_2\neq k_1 } M^2 &\lesssim  \sup_{k,\tau}  \sum_{k_1\neq k,k_2\neq k_1 } \frac{\la k \ra^{2s+2a} \la k_1\ra^{-2s } \la k_2\ra^{-2s } \la k-k_1+k_2\ra^{-2s }  }{  \la \tau- k^2\ra^{2-2b} \la \tau-k_1^2+k_2^2-(k-k_1+k_2)^2\ra^{4b-1}}\\
&\lesssim \sup_{k }  \sum_{k_1\neq k,k_2\neq k_1 }\frac{\la k \ra^{2s+2a} \la k_1\ra^{-2s } \la k_2\ra^{-2s } \la k-k_1+k_2\ra^{-2s }  }{   \la k^2-k_1^2+k_2^2-(k-k_1+k_2)^2\ra^{2-2b}}.
\end{align*}
The last line follows by    the simple fact
\be\label{trivial}
\la \tau-n\ra \la \tau-m\ra \gtrsim \la n-m\ra.
\ee
Since we have only the nonresonant terms, it suffices to estimate
$$
 \sum_{k_1, k_2  }\frac{\la k \ra^{2s+2a} \la k_1\ra^{-2s } \la k_2\ra^{-2s } \la k-k_1+k_2\ra^{-2s }  }{   \la k-k_1\ra^{2-2b} \la  k_1-k_2 \ra^{2-2b}}.
$$
To estimate this sum we consider the cases $|k-k_1+k_2|\gtrsim |k|$,  $|k_1|\gtrsim |k|$, and $|k_2|\gtrsim|k|$. In the estimate we always take $0<b-\frac{1}{2}$ sufficiently small for any given $s>0$ and $0\leq a<\min(2s,\frac{1}{2})$.

In the first case, using Lemma~\ref{lem:sums} we bound the sum by
$$
\sum_{k_1,k_2}\frac{\la k \ra^{ 2a} \la k_1\ra^{-2s } \la k_2\ra^{-2s }    }{   \la k-k_1\ra^{2-2b} \la  k_1-k_2 \ra^{2-2b}}\les \sum_{k_1}\frac{\la k \ra^{ 2a}   \phi_{\max(2s,2-2b)}(k_1)  }{   \la k-k_1\ra^{2-2b} \la  k_1  \ra^{2s+\min(2-2b,2s)}}.
$$
In the case $s\geq \frac{1}{2}$, we bound the sum by
$$
 \sum_{k_1}\frac{\la k \ra^{ 2a}   \log(1+\la k_1\ra) }{   \la k-k_1\ra^{2-2b} \la  k_1  \ra^{2s+2-2b }}\les \la k\ra^{2a-2+2b+}\les 1
$$
provided $a<\frac{1}{2}$.

In the case $0<s<\frac{1}{2}$, using Lemma~\ref{lem:sums}    we bound the sum by
$$
\sum_{k_1}\frac{\la k \ra^{ 2a}     }{   \la k-k_1\ra^{2-2b} \la  k_1  \ra^{4s+ 1-2b }}\les \left\{\begin{array}{lc}\la k\ra^{2a+4b-4s-2} & 0<s\leq \frac{1}{4},\\
                                                                                                    \la k\ra^{2a+2b -2} & \frac{1}{4}<s < \frac{1}{2}.
                                                                                                   \end{array}\right.
$$
This is bounded in $k$ provided that $0\leq a<\min(2s,\frac{1}{2})$.

The second case is identical to the first case after renaming the variables: $n_1=k-k_1+k_2$, $n_2=k_2$.

In the third case, after renaming the variables $n_1=k_1$, $n_2=k-k_1+k_2$, and using Lemma~\ref{lem:sums} we bound the sum by
$$
\sum_{ n_1,n_2 }\frac{\la k \ra^{ 2a }  \la n_1\ra^{-2s} \la n_2\ra^{-2s}   }{   \la k-n_1\ra^{2-2b} \la k-n_2\ra^{2-2b}  }\les \la  k\ra^{2a}  \phi^2_{\max(2s,2-2b)}(k)  \la  k   \ra^{-2\min(2-2b,2s)}\les 1,
$$
provided that  $s>0$ and $0\leq a<\min(2s,1)$.
This finishes the proof of Proposition~\ref{lem:smooth}.

\vskip 0.1in
\noindent
{\bf Remark 2.} A similar smoothing property was established in \cite{erdtzi} for the KdV equation. However, in that case, smoothing is not immediate within the context of the $X^{s,b}$ norms as the following example shows \cite{kpv}. Fix $s\in \R$. Consider $\widehat u(n,\tau)= \delta(n-M) \chi_{[-1,1]}(\tau-n^3)$, and
$\widehat v(n,\tau)=\delta(n-1) \chi_{[-1,1]}(\tau-n^3)$ with $M\gg 1$. The inequality
$$
\|(uv)_x\|_{X^{s+a,-1/2}_{KdV}}\lesssim \|u\|_{X^{s,1/2}_{KdV}} \|v\|_{X^{s,1/2}_{KdV}}
$$
fails for any $a>0$.

The analogue of Theorem~\ref{mainpro} was proved in \cite{erdtzi} after transforming the KdV equation using differentiation by parts. However, the dimension statement in Theorem~\ref{mainthm} is open for KdV since the corresponding result for the linear part (Airy equation) is not known. We will address this issue in future work.

\vskip 0.1in
\noindent
{\bf Remark 3.} The gain of almost half a derivative in Theorem~\ref{mainpro} can be inferred by taking the first Picard iteration for the cubic NLS and test the term
$$
e^{it\partial_{xx}}\int_0^t e^{ -it^\prime \partial_{xx}} \big[e^{it^\prime \partial_{xx}} g \overline{\big(e^{it^\prime \partial_{xx}} g\big)}e^{it^\prime \partial_{xx}} g\big] dt^\prime.
$$
On the Fourier side, ignoring  resonant terms, one sees that for $g\in L^2$ the term
$$\sum_{k_1+k_2+k_3=k} \int_0^t e^{-2i(k_1+k_2)(k_2+k_3)t^\prime} \widehat{g}(k_1)\widehat{g}(k_2)\widehat{g}(k_3)dt^\prime$$
$$= \sum_{k_1+k_2+k_3=k}\frac{\widehat{g}(k_1)\widehat{g}(k_2)\widehat{g}(k_3)}{2(k_1+k_2)(k_2+k_3)} (e^{-2i(k_1+k_2)(k_2+k_3)t}-1)
$$
is in $H^{\frac{1}{2}-}$ but not necessarily better.

\end{document}